\input amssym.def
\input amssym
\magnification=\magstep1
\def\nbigskip{\bigskip\noindent}
\def\bbigskip{\bigskip \bigskip  }
\def\nbbigskip{\bbigskip \noindent }
\def\nbbbigskip{\bigskip\bigskip\bigskip\noindent}
\def\nmedskip{\medskip\noindent}
\def\qq{/\kern-.3em /}
\def\C{{\Bbb C}}
\def\R{{\Bbb R}}
\def\im{{\rm Im\,}}
\def\Sl{{\rm SL}}
\def\Su{{\rm SU}}
\def\buildunder#1#2{\mathrel{\mathop{\kern0pt #2}\limits_{#1}}}

\def\qed{\hfill$\square$}
\def\int{{\rm int}}

\def\H{{\bf H}}
\def\tube{{\cal T}}
\def\litem{\par\noindent\hangindent=\parindent\ltextindent}
\def\ltextindent#1{\hbox to \hangindent{#1\hss}\ignorespaces}
\def\litem{\par\noindent\hangindent=\parindent\ltextindent}
\def\ltextindent#1{\hbox to \hangindent{#1\hss}\ignorespaces}

\centerline{\bf The minimum principle from a Hamiltonian point of view}

\nbigskip
\centerline{\sl Peter Heinzner}

\nbbbigskip
This paper grew out of a manuscript written by Xiang-Yu Zhou
on the extended future tube conjecture ([Z]). There the following 
problem is considered. 
Let $<\ ,\ >$ denote the Lorentz product on $\R^4$ and also its 
$\C$-bilinear extension to $\C^4$. The future tube 
$\tube$ is by definition the tube domain in $\C^4=\R^4+i\R^4$ over
the positive light cone $C^+=\{y=(y_0, y_1,y_2,y_3)\in \R^4;\ y_0>0,\ 
<y,y>=(y_0)^2-(y_1)^2-(y_2)^2-(y_3)^2 > 0\}$, i.e.,
$$\tube=\{z \in \C^4;\ \im z \in C^+\}\,.$$ 
This domain is invariant under the action of the connected component
$G_\R$ of the identity of the homogeneous Lorentz group 
${\rm O}_\R(1,3)$. Now consider  
the $N$-fold product $\tube^N$ with the diagonal action of    
$G_\R$. The extended future tube $(\tube^N)^\C$ is by definition
the orbit of $\tube^N$ under the action of the complexified group 
$G$ of $G_\R$.  In other words 
$$(\tube^N)^\C=
G\cdot \tube^N=\{(g\cdot z_1,\ldots, g\cdot z_N);\ g\in G,\ z_j\in \tube\}\,. $$              Note that $G$ is the group ${\rm SO}_4(\C)$ which is 
defined by the quadratic form $<\ ,\ >$.      

In section \S3 we give a conceptual proof of the extended future tube conjecture.

\nbigskip
{\bf Theorem 1.} {\it The extended future tube is  
a domain of holomorphy.}

\nbigskip
This result has conjecturally been known in quantum field theory for more then
thirty years. For its relevance and other publications concerning 
problems related to it we refer the reader to the literature
([B-L-T], [H-S], [J], [S-W], [S-V]).

\smallskip
Let us put the context of Theorem 1 in a more general framework
which has been studied intensively over the last years from the
point of view of holomorphic group actions.   

\smallskip
Let $G_\R$ be a connected real form of a  complex Lie group 
$G$ and let $X$ be a $G_\R$-stable domain of holomorphy  in $Z$ 
such that $G\cdot X=Z$. 
Under which conditions on $X$ is $Z$ the natural domain of definition 
of the $G_\R$-invariant holomorphic functions on $X$? 

\smallskip
If $Z$ is a submanifold of a Stein manifold
and every $G_\R$-invariant holomorphic function  extends
to $Z$, then it also extends to the envelope of holomorphy of $Z$. Thus
one also has to ask under which additional requirements is $Z$ a Stein 
manifold. 

\smallskip
In order that an invariant holomorphic function extends to $Z=G\cdot X$ it 
is sufficient that $X$ is orbit connected, i.e., for every $z\in Z$ the 
set $\{g\in G;\ g\cdot z\in X\}$ is connected (see [H]). This condition is
satisfied for example for the extended future tube ([H-W]). Thus under this 
condition the question is whether  $Z$ is a Stein manifold.

\smallskip  
If $Z$ is a domain in a Stein manifold $V$, then $Z$ itself 
is a Stein manifold 
if one can find a plurisubharmonic function $\Psi$ on $Z$ which goes to
$+\infty$ at every boundary point of $\partial Z\subset V$. 
There is a very natural way to construct $G$-invariant plurisubharmonic   
functions out of $G_\R$-invariant functions on $X$ which was first proposed 
by Loeb in his fundamental paper [L]. In this paper Loeb used an
extended version of Kiselman's minimum principle in order to construct
invariant plurisubharmonic function. The main idea is the following. 
Assume that there is a nice quotient $\pi:Z\to Z/G$ and let $\phi$ be a 
smooth $G_\R$-invariant plurisubharmonic function on $X$ which is a strictly
plurisubharmonic  exhaustion on each fibre of $\pi\vert X$. Then the 
fibre wise minimum 
of $\phi$ defines a function $\psi$ on $Z/G$ which is a candidate for
a plurisubharmonic function on $Z/G$. This procedure can be 
described in terms of Hamiltonian actions as follows. 

\smallskip
Assume for simplicity that $\phi$ is strictly plurisubharmonic. Then 
$\omega:=2i\partial\bar\partial\phi$ defines an invariant K\"aher form
on $X$ and $\mu(x)(\xi)=d\phi(J\xi_X)$ is the associated 
moment map $\mu:X\to \frak g^*_\R$. 
In this situation $\mu^{-1}(0)$ is the set of fibre wise critical points of 
$\phi$ which in good cases are exactly the points such that the restriction 
of $\phi$ to the fibre attains its minimum. Again under some 
additional assumption, it then follows from the principle of 
symplectic reduction that the reduced 
space $\mu^{-1}(0)/G_\R$ has  a symplectic structure which in fact is 
K\"ahlerian and moreover is given by the function 
$\psi$ which is induced on $\mu^{-1}(0)/G_\R$
by $\phi\vert \mu^{-1}(0)$. It turns out that in 
the situation under consideration
the procedures given by symplectic reduction 
and minimum principle are compatible. This is well known in the case where 
$G_\R$ is a compact Lie group (see [H-H-L], where a much more general result is
proved) and we give here precise  conditions such that it also works for a 
non compact group $G_\R$.     

\smallskip
Although there is no geometric quotient of $Z$ in the case of 
the extended future tube, we have a quotient 
$\pi:(\C^4)^N\to (\C^4)^N\qq G$ which is given by the 
invariant holomorphic functions on $(\C^4)^N$  and it is a fundamental 
fact that the 
extended tube $(\tube^N)^\C$ is saturated with respect to $\pi$ ([H-W],
see \S3 for additional remarks).
In this case it turns out that this invariant theoretical quotient has 
sufficiently many good properties in order to apply the main result 
of this paper which we formulate now.

\smallskip
Let $V$ be a holomorphic Stein $G$-manifold such that there exists 
almost a  quotient $\pi:V\to V\qq G$. More precisely we will assume 
that $V\qq G$ is a complex space, $\pi:V\to V\qq G$ is a 
$G$-invariant surjective holomorphic map and for an analytically 
Zariski open $\pi$-saturated subset $V^0$ of $V$ the restriction map
$\pi:V^0\to V^0\qq G$ is a holomorphic fibre bundle with typical fibre 
$G/H$. Let $X$ be a $G_\R$-stable domain in $V$ such that 
$Z:=G\cdot X$ is saturated with respect to $\pi:V\to V\qq G$.

\nbbigskip
{\bf Theorem 2.} {\it Let  $\phi:X\to \R$ be a smooth non-negative 
$G_\R$-invariant plurisubharmonic function and assume that 

{\parindent=1cm
\smallskip
\litem{{\rm(i)}} $X^0:=X\cap V^0$ is orbit connected, $X^0:=X\cap V^0$, 

\smallskip
\litem{{\rm(ii)}} the restriction of $\phi^0:=\phi\vert X^0$ to the part of 
the $\pi$-fibres which lies in $X^0$ is strictly plurisubharmonic, 

\smallskip
\litem{{\rm(iii)}} $\phi^0$ is proper  mod $G_\R$ along $\pi\vert Z^0$ where
$Z^0:=V^0\cap Z$ and 
\smallskip
\litem{{\rm(iv)}} $\phi$ is a weak exhaustion of $X$ over $V\qq G$,
}

\smallskip\noindent
Then $Z=G\cdot X$ is a Stein manifold.}

\nbbigskip
In the case where $G_\R$ acts properly on $X^0$ condition (iii) 
means that the map $\phi^0\times \pi\vert X^0:X^0\to 
\R\times (Z^0\qq G)$ induces a proper map $X^0/G_\R\to\R\times (Z^0\qq G)$.
By a weak exhaustion of $X$ over $V\qq G$ we mean a function which goes
to $+\infty$ on a sequence if the corresponding sequence in $V\qq G$ 
converges to a boundary point of $Z\qq G$ in $V\qq G$.  

\smallskip
In the case where the $G$-action on $Z^0$ is assumed to be free, 
the theorem can be proved rather directly by applying Loeb's 
minimum principle. For a compact group it is a consequence of the
methods presented in [H-H-K] (see also [H-H-L]).

\smallskip
In the last section we recall a result in [Z] about the 
orbit geometry of the extended future tube together with some
previously known facts in order to verify that the conditions of  
Theorem 2 are satisfied in the case of the extended future tube.
Orbit connectedness as well as saturatedness of the extended tube
is proved in [H-W]. 

\smallskip
In the meantime Zhou informed me that he realized that the 
properness property of $\phi$ is one of the crucial points in the
application of the minimum principle and that he has taken this
into account in a revised version of his manuscript. 
The last step in the proof of the extended future tube conjecture 
should therefore be assigned to him.

\nbbbigskip
{\bf 1. Hamiltonian actions on K\"ahler spaces.}

\nbbigskip
Let $G_\R$ be a real  connected Lie group and 
$X$ a complex $G_\R$-space, i.e., $G_\R$ acts on $X$ 
by holomorphic transformations such that the action 
$G_\R\times X\to X, \ (g,x)\to g\cdot x$, 
is real analytic. If $\omega$ is a smooth $G_\R$-invariant K\"ahler
structure  on $X$, then a $G_\R$-equivariant smooth map $\mu$ 
from $X$ into the dual $\frak g_\R^*$ of the Lie
algebra $\frak g_\R$ of $G_\R$ is said to be an equivariant moment map
if 
$$d\mu_\xi=\imath_{\xi_X}\omega$$
holds on every $G_\R$-stable complex 
submanifold $Y$ of $X$. Here $\omega$ denotes the K\"ahler form 
on $Y$ induced by the K\"ahlerian structure on $X$ (see [H-H-L]), 
$\mu_\xi:=<\mu,\xi>$ is the component of $\mu$ in 
the direction of $\xi\in \frak g_\R$, $\xi_X$ is the vector field on $X$ induced 
by $\xi$ and  $\iota_{\xi_X}\omega$ denotes the one form given by contraction,
i.e., $\eta\to\omega(\xi,\eta)$.

\nmedskip
{\it Example.} If $\omega$ is given by a smooth strictly plurisubharmonic
$G_\R$-invariant function $\phi$, i.e., $\omega=2i\partial \bar\partial \phi$ 
on every smooth part of $X$, then 
 
$$\mu_\xi(x):=d\phi(J\xi_X)=(i(\partial-\bar\partial)\phi)(\xi_X)=d^c\phi(\xi_X)$$
defines an equivariant moment map. This follows from invariance of $\phi$, since
in this case we have 

$$d\mu_\xi=d\imath_{\xi_X} d^c\phi=-\imath_{\xi_X} d d^c\phi = 
\imath_{\xi_X}2i\partial\bar\partial\phi\,.$$
Here we use the formula 

$${\cal L}_\xi\alpha = \imath_\xi d \alpha+d\imath_\xi \alpha$$ 
for all vector fields $\xi$ and differential forms $\alpha$ where
${\cal L}_\xi$ denotes the Lie derivative in the direction of
$\xi$.

\medskip
Later we will need the following fact about the zero level set of $\mu$.

\nbbigskip
{\bf Lemma.} {\it Assume that $X$ is smooth and that $G_\R$ acts
properly on $X$.  If the dimension of the $G_\R$-orbits 
in $\mu^{-1}(0)$ is constant, then $\mu^{-1}(0)$ is a 
submanifold of $X$.}

\nbigskip
{\it Proof.} Since the action is assumed to be proper, there is 
a local normal form for the moment map (see e.g. [A] or [H-L]). 
The statement is an easy consequence of this fact (see e.g. [A]. In
[S-L] the argument is given for a compact group $G_\R$). 
\qed

\nbbigskip
{\it Remark$\,$1.} It can be shown that the converse of the Lemma 
also holds. We will not use this fact here.

\nmedskip
{\it Remark$\,$2.} The properness assumption is very often 
satisfied. Since one may assume that $G_\R$ acts effectively,
$G_\R$ is a Lie subgroup 
of the group $I$ of isometries of the Riemannian manifold $X$. 
The group of isometries acts properly on
$X$ and consequently  the $G_\R$-action on $X$ is proper if 
and only if $G_\R$ is a closed
subgroup of $I$. This is the case if and only if there is a point
$x\in X$ such that $G_\R\cdot x$ is closed and the isotropy group 
$(G_\R)_x:=\{g\in G_\R;\, g\cdot x=x\}$ is compact.

\nmedskip
{\it Remark 3.} If $G_\R$  acts such that the isotropy groups are 
discrete, then $\mu$ has maximal rank. Thus in this case $\mu^{-1}(0)$ 
is obviously a submanifold of $X$. Moreover $T_x(\mu^{-1}(0))=\ker d\mu(x)$
for all $x\in \mu^{-1}(0)$.

\nbbbigskip
{\bf 2. Hamiltonian actions on invariant domains}

\nbbigskip
Let $G$ be a connected complex Lie group and $Z$ a holomorphic $G$-space, i.e.,
the action $G\times Z\to Z$ is assumed to be a holomorphic map.
Let $G_\R$ be a connected real form of $G$. By an invariant domain in
$Z$ we mean in the following a $G_\R$-stable connected open
subspace $X$ of $Z$. In the homogeneous case we have the following

\nbbigskip
{\bf Lemma 1.} {\it Let $X$ be an invariant domain in $Z$ and assume 
that $Z$ is $G$-homogeneous. If the zero level set of 
$\mu:X\to \frak g_\R$ is  
not empty, then $\mu^{-1}(0)$ is a Lagrangian submanifold of $X$ and 
each connected component of $\mu^{-1}(0)$ is a $G_\R$-orbit.}

\nbigskip
{\it Proof.} For $z_0\in X$ let $N$ be an open convex neighborhood of 
$0\in \frak g_\R$ such that $U:=G_\R \cdot \exp iN \cdot z_0 \subset X$. 
Since $G_\R\cdot \exp iN$ is a neighborhood of $G_\R$ in $G$, the set 
$U$ is a neighborhood of $G_\R\cdot z_0$ in $X$.
The proof of Lemma 1 is a consequence of the following

\nmedskip
{\it Claim.} $U\cap \mu^{-1}(0)=G_\R\cdot z_0$ for $z_0\in \mu^{-1}(0)$.

\nmedskip
In order to proof the claim, let $z\in U\cap \mu^{-1}(0)$ be given.
Then there are $h\in G_\R$ and $\xi\in N$ such that 
$z=h\exp i\xi\cdot z_0\in U\cap  \mu^{-1}(0)$. Thus $\exp i\xi\cdot z_0\in
\mu^{-1}(0)\cap U$ and $z_t:=\exp it\xi\cdot z_0\in U$ for $t\in [0,1]$. Note
that $J\xi_X(x)={d\over dt}\vert_{t=0}\exp it\xi\cdot x$ is the gradient 
flow of $\mu_\xi$ with respect to the Riemannian metric induced by $\omega$.
Thus, if $z_t$ is not constant, then $t\to \mu_\xi(z_t)$ is strictly 
increasing. This contradicts $\mu_\xi(z_0)=0=\mu_\xi(z_1)$. Therefore 
$z_0=\exp it\xi\cdot z_0$ 
for all $t\in \R$. This implies $z=h\cdot z_1=h\cdot z_0\in G_\R\cdot z_0$.

\smallskip
It is a consequence of the claim that every $G_\R$-orbit is closed in $X$.
Therefore every component of $\mu^{-1}(0)$ is a $G_\R$-orbit.
It remains to show that these orbits are Lagrangian. 
Since $\mu(G_\R\cdot z_0)=0$ we have 
$$0=d\mu_\xi(\eta_X(z_0))=\omega(\xi_X (z_0),\eta_X (z_0))$$ 
for all $\xi,\eta\in \frak g_\R$. This 
means that $G_\R\cdot z_0$ is an isotropic submanifold of $X$. 
In particular, $\dim_\R\, G_\R\cdot z_0\le \dim_\C \,X$. 
In general the tangent space $T_{z_0}(G_\R\cdot z_0)$ spans 
$T_{z_0}X$ over $\C$. Thus 
$\dim_\R G_\R\cdot z_0\ge \dim_\C G\cdot z_0=\dim_\C X$. This
shows that $\dim_\R G_\R\cdot z_0={1\over 2}\dim_\R X$. 
Hence $G_\R\cdot z_0$ is Lagrangian.
\qed

\bbigskip
Every Lagrangian submanifold of a K\"ahler manifold 
is totally real. Thus, if $Z$ is $G$-homogeneous, then 
$\mu^{-1}(0)$ is a totally real submanifold of $X$. Note that
the $G_\R$-orbits in $\mu^{-1}(0)$ are closed since they are
connected components of the zero fibre of $\mu$. 
Now if $G_\R$ is such that $0\in \frak g_\R^*$ is the only $G_\R$-fixed 
point, then $x\in \mu^{-1}(0)$ if and only if the orbit $G_\R\cdot x$ 
is isotropic. This condition holds for example for a semisimple Lie
group. 

\smallskip
It almost never happens that there is a $G_\R$-invariant K\"ahler form
$\omega$ which is defined on $Z$. For example, if $G_\R$ is a simple non 
compact Lie group or more generally a semisimple Lie group without 
compact factors, then there does not exist a $G_\R$-invariant K\"ahler 
form on a non trivial holomorphic $G$-manifold $Z$. In order to see this, recall
that since $G_\R$ is semisimple there is a moment map $\mu:Z\to \frak g^*_\R$.
Now let $\frak g_\R= \frak k \oplus \frak p$ be a Cartan decomposition where 
$\frak k$ is the Lie algebra of the maximal compact subgroup of $G_\R$. Then
$\frak u=\frak k \oplus i\frak p$ is the Lie algebra of the maximal 
compact subgroup $U$
of $G$. For $\xi\in i \frak p$ the image of the one-parameter group 
$\gamma:t\to \exp it\xi$ lies
in $U$ and therefore there is a basis of $\frak p$ consisting of $\xi$'s such that 
the image of $\gamma$ is compact, i.e., isomorphic
to $S^1$. But $\gamma$ is the flow of the 
gradient vector field of $\mu_\xi$ and therefore 
$t\to \mu_\xi(\gamma(t)\cdot z)$ is strictly
increasing for every $z\in Z$. This implies that $\gamma$ acts trivially on 
$Z$. Since $G$ is semisimple 
and contains no compact factor, $G$ itself is the 
smallest complex
subgroup of $G$ which contains $\exp \frak p$. Thus $G$ acts trivially on $Z$.

\smallskip
A geometric  interpretation of the zero fibre $\mu^{-1}(0)$  of an 
equivariant moment map $\mu:X\to \frak g_\R^*$ associated to a 
smooth $G_\R$-invariant strictly plurisubharmonic function 
$\phi:X\to \R$ (see Section 1, Example) can be given in the
case where $X$ is an invariant domain in $Z$ as follows. For 
$x\in X$ let $\Omega(x):=\{g\cdot x;\, g\in G \ \hbox{and}\ g\cdot x \in X\}$
be the local $G$-orbit of $G$ through $x$ in $X$ where 
$(g,x)\to g\cdot x$ denotes the $G$-action on $Z$. Then by 
$G_\R$-invariance of $\phi$ we have $$\mu^{-1}(0)=\{x\in X;\,  
\ x\ \hbox{is a critical point of}\ \phi\vert \Omega(x)\}\,.$$
 
\smallskip 
We consider now invariant domains $X$ in $G$-homogeneous spaces $Z$ such that
there is a moment map associated to $\phi:X\to \R$ more closely. In order to 
do that we first introduce the notion of an exhaustion mod $G_\R$.

\smallskip 
Let $F$ be a complex space with a proper $G_\R$-action and let $F/G_\R$ 
be the space of $G_\R$-orbits endowed with the quotient topology. 
A $G_\R$-invariant function $f:F\to \R$ is said to be proper 
mod $G_\R$ if the induced map
$\bar f:F/G_\R\to \R$ is proper. The map $f$ is said to be 
an exhaustion mod $G_\R$ if $\bar f$ is an exhaustion, i.e., if
for all $r\in \R$ the set $\{q\in F/G_\R;\, \bar f(q)< r\}$ is
relatively compact in $F$. Note that a $G_\R$-invariant continuous function 
which is bounded from below is proper mod $G_\R$ if and only if it 
is an exhaustion mod $G_\R$.

\nbbigskip
{\bf Lemma 2.} {\it Let $Z$ be $G$-homogeneous and assume that the 
$G_\R$-action on $X$ is proper. Let 
$\phi:X\to \R$ be a smooth strictly plurisubharmonic $G_\R$-invariant 
function which is an exhaustion mod $G_\R$. Then there is a $z_0\in X$
such that 
$$G_\R\cdot z_0=\mu^{-1}(0)=\{z\in X;\ 
\phi(z) \ \hbox{is a minimal value of} \ \phi\}\,.$$}

\nbigskip    
{\it Proof.}
Since $\phi$ is plurisubharmonic and an exhaustion mod $G_\R$
there is a point $z_0\in X$ which is a minimum for $\phi$. In particular,
$\mu^{-1}(0)$ is not empty where $\mu$ denotes the moment map associated
with $\phi$.  We have to prove that $\mu^{-1}(0)$ is connected.
By Lemma 1, every  connected component of the
set $M_\phi=\mu^{-1}(0)$ of critical points of $\phi$ is a  $G_\R$-orbit. 
We claim that the $G_\R$-orbits
are non degenerate in the sense that the Hessian of $\phi$ in 
normal directions is positive definite. This is seen as follows.

\smallskip
The vector fields $J\xi_X$, $\xi \in \frak g_\R$ span 
the normal space at $x\in M_\phi$
and
$$(J\xi_X)(J\xi_X(\phi))=\imath_{J\xi_X}d\mu_\xi=\omega(\xi_X,J\xi_X)\,.$$
Hence the Hessian at $x\in M_\phi$ is positive in the normal directions.
Since $\phi$ is proper mod $G_\R$ and the gradient vector field of $\phi$ with
respect to the $G_\R$-invariant K\"ahler metric given by 
$2i\partial\bar\partial\phi$
is $G_\R$-invariant, Lemma 2 follows from standard 
arguments in Morse Theory.
\qed

\nbbigskip
In the situation of Lemma 2 
every critical point of $\phi$ is a minimum and the
set of these points is a $G_\R$-orbit  and coincides with 
$\mu^{-1}(0)$. 
 
\smallskip
We will now generalise the results in the homogeneous case to spaces 
$Z$ which possess a geometric $G$-quotient and $X$ is a weakly orbit
connected invariant domain in $Z$. 
Here a $G_\R$-stable subset $X$ of $Z$ is
said to be weakly orbit connected if 
for every $x\in X$ the local $G$-orbit $\Omega(x):=\{g\cdot x\in X;\, g\in G\}$
is connected.

\nmedskip
{\it Remark 1.} A $G_\R$-invariant set $X$ in $Z$ is said to be orbit connected
if for every $x\in X$ the set $\Omega_x:=\{g\in G;\, g\cdot x\in X\}$ is
connected. This is a stronger concept then weakly orbit connectedness.

\medskip
Let $Z$ be a holomorphic $G$-space such that there is a geometric quotient
$\pi:Z\to Z/ G$. By this we mean that the orbit space $Z/G$ is a complex space
such that the quotient map $\pi:Z\to Z/G$ is holomorphic. Moreover we assume
that the structure sheaf of $Z/G$ is the sheaf of invariants, i.e., 
for an open subset $Q$ of $Z/G$ a function
$f:Q\to \C$ is holomorphic if and only if $f\circ \pi:\pi^{-1}(Q)\to \C$ is 
holomorphic. 

\smallskip
Now let $X\subset Z$ be an invariant domain which lies surjectively over 
$Z/ G$ or 
equivalently such that $Z=G\cdot X$. Assume that $G_\R$ acts properly on 
$X$ and that $X$ is weakly orbit connected. Let $\phi:X\to \R$ be a 
smooth $G_\R$-invariant strictly plurisubharmonic function which is 
an exhaustion mod $G_\R$ along $\pi$, i.e., $\pi^{-1}(C)\cap \{x\in X;\
\phi(x)\le r\}/G_\R\subset X/G_\R$ is compact for every compact subset 
$C$ in $Z/G$ and $r\in\R$. We set $M_\phi=\mu^{-1}(0)$ where 
$\mu:X\to \frak g_\R$ denotes the moment map associated with $\phi$.

\nbbigskip
{\bf Proposition 1.} {\it The map $\bar \imath :M_\phi/G_\R\to Z/ G$ induced
by the inclusion $\imath:M_\phi\to Z$ is a homeomorphism. If $X$ is  
a manifold, then $M_\phi$ is smooth and 
$$T_xM_\phi=\ker d\mu(x)$$  
holds for all $x\in M_\phi$.}

\nbigskip
{\it Proof.} The map $\bar \imath$ is continuous and by Lemma 2 it 
is also a bijection. We claim that $\bar\imath$ is proper. Since the
$G_\R$-action on $M_\phi$ is proper, $M_\phi/G_\R$ is a locally compact
topological space. Thus properness of $\bar \imath$ implies that 
$\bar\imath$ is a homeomorphism.

\smallskip
Let $(q_n)$ be a sequence in $M_\phi/G_\R$ and $x_n$ a point in
$M_\phi$ which lies over $q_n$. Assume that $(\pi(x_n))=(\bar\imath(q_n))$
has a limit in $Z/G$ and let $x_0\in M_\phi$ be a point which lies over
$\lim \pi(x_n)$. If some subsequence of $\phi(x_n)$ goes to infinity, 
then we may assume $\phi(x_n)>\phi(x_0) + 1$ for all $n$. 
Since $\pi:Z\to Z/G$ is an open map, there are $g_n\in G$ such that 
$\lim g_n\cdot x_n=x_0$ for some subsequence. This is a contradiction since 
$\phi(x_n)< \phi (g_n\cdot x_n)$ for all $n$ such that $g_n\cdot x_n\in X$.
Thus, since $\phi$ is assumed to be an exhaustion mod $G_\R$ 
along $\pi$,
there are $h_n\in G_\R$ such that a subsequence of
$(h_n\cdot x_n)$  converges to $x_0$. This implies that a subsequence of
$(q_n)$ converges in $M_\phi$. So far we proved that $\bar \imath$ is a 
homeomorphism.

\smallskip
Assume now that $X$ is smooth. The existence of a geometric quotient implies
that the dimension of the $G$-orbits in $Z$ is constant and therefore 
this is also true for the $G_\R$-orbits in $M_\phi$ (Lemma 1). 
Thus $M_\phi$ is a submanifold of $X$ (Section 1, Lemma). 
Since $T_x M_\phi$ is a subspace of $ \ker d\mu(x)$ 
and $\ker d\mu(x)=T_x (G_\R\cdot x) + T_x (G\cdot x)^\perp$, the claim follows
from the obvious dimension count as follows. Let $d:=\dim_\R G_\R\cdot x$
for $x\in M_\phi$. Note that $d$ is the complex dimension of the 
$\pi$-fibres. Thus
$\dim_\R M_\phi = \dim_\R M_\phi/G_\R + d = \dim_\R Z / G + d = 
\dim_\R T_x(G\cdot x)^\perp + \dim_\R G_\R\cdot x$
implies that $T_x M_\phi = \ker d\mu (x)$ for all $x\in M_\phi$.
\qed

\nbbigskip
{\it Remark 2.} Without a reference to an embedding into a  
holomorphic $G$-space one can show that $\mu^{-1}(0)/G_\R$ 
is a complex space in a natural way (see [A-H-H] and [A]).

\medskip
If $G_\R$ does not act properly on $X$, then let
$\overline {G}_\R$ be the closure of $G_\R$ in the group $I$ of
isometries of the K\"ahler manifold $X$. Since the $G_\R$-orbits
in $M_\phi=\mu^{-1}(0)$ are closed (Lemma 1), it follows that they coincide 
with the $\overline {G}_\R$-orbits. Moreover 
$\phi$ is $\overline {G}_\R $-invariant
and $M_\phi=\bar\mu^{-1}(0)=:\overline M_\phi$, where $\bar \mu$ is the moment map 
associated with $\phi$. 
Now if one redefines an exhaustion mod $G_\R$ along $\pi$ in terms of
sequences in $X$, then also in this case $M_\phi$ is smooth and
$T_xM_\phi=T_x(G_\R\cdot x)\oplus T_x(G\cdot x)^\perp =\ker d\mu(x)$  
holds for all $x\in M_\phi$.

\smallskip
Proposition 1 can be 
generalised to the case where $\phi:X\to \R$ is only assumed to be 
plurisubharmonic and strictly plurisubharmonic on the
fibres. More precisely we have the following consequence  which
can be thought of as a version of Loeb's minimum
principle (see [L]).

\nbbigskip
{\bf Corollary 1.} {\it Let $X\subset Z$ be a weakly orbit connected 
invariant domain with $\pi(X)=Z$ and $\phi:X\to \R$ a smooth $G_\R$-invariant 
plurisubharmonic function which is an exhaustion mod $G_\R$ along 
$\pi$ such that the restriction of $\phi$ to the local 
$G$-orbits in $X$ is a strictly plurisubharmonic exhaustion mod
$G_\R$. If $\pi:Z\to Z/G$ is a holomorphic bundle, then

{\parindent=1cm
\smallskip
\litem{{\rm(i)}} $M_\phi=\mu^{-1}(0)$ is smooth where $\mu:X\to \frak g_\R^*,\,
\mu_\xi=d\phi(J\xi_X)$,

\smallskip
\litem{{\rm (ii)}}  $T_x(M_\phi)=\ker d\mu(x)$ for all $x\in M_\phi$.

\smallskip
\litem{{\rm(iii)}} $M_\phi/\overline{G}_\R $ is homeomorphic to $Z/G$ and
the function $\psi:Z/G\to \R$ which is induced by $\phi\vert M_\phi$
is a smooth plurisubharmonic function.\par}
}

\nbigskip
{\it Proof.} We may assume that $G_\R$ acts properly on $X$ and, since the 
statements are local over $Z/G$ that $Z/G$
is a Stein manifold. Let $\rho:Z\to \R$ be the the pull back of
a strictly plurisubharmonic function on $Z/G$. Then $\phi + \rho$ is 
$G_\R$-invariant, strictly plurisubharmonic and an exhaustion mod $G_\R$ 
on the local $G$-orbits in $X$. Since $d\rho(J\xi_X)=0$ for all $\xi \in \frak g_\R$,
the moment map associated with $\phi + \rho $ is the same as the 
moment map associated with $\phi$. Thus Proposition 1 implies directly
(i), (ii) and the first part of (iii). It remains to show that 
$\psi:Z/G\to \R$ is a smooth plurisubharmonic function. 

\smallskip
For the plurisubharmonicity of $\psi$ we recall 
the calculation in [H-H-L], \S2. For $z\in M_\phi$ we have
$T_z(M_\phi)=\ker d\mu(z)=T_z(G_\R\cdot z)\oplus T_z(G\cdot z)^\perp$. 
We may assume that $Z=G/H\times \Delta$ where $\Delta$ 
is an open neighborhood of $0$ in $\C^d\cong T_z(G\cdot z)^\perp$, 
and $\pi(z)=0$  where $\pi$ is given by the projection on the second factor.
Furthermore there is a section 
$\eta:\Delta\to M_\phi,\, \eta(w)=(\sigma(w),w)$ 
and therefore we have $\psi(w)=\phi(\eta(w))$. A direct calculation 
shows that 
$$\partial\bar\partial\,\psi(0)=\partial\bar\partial\,\phi(\eta(0))\,.$$
Here one has to use that $d\phi(z)=0$ and that $d\sigma(0)=0$.
Thus $\psi$ is plurisubharmonic and smooth.
\qed

\bbigskip
If $\phi$ is strictly plurisubharmonic, then the proof shows 
that $\psi$ is also strictly plurisubharmonic. For a proper $G_\R$-action
the space $Z/G$ is
then given by symplectic reduction $M_\phi/G_\R$ and  the induced K\"ahlerian
structure on $Z/G$ is determined by the function 
$\psi(q)=\buildunder{x\in \pi^{-1}(q)\cap X}\inf \phi(x)$ 
which is obtained by applying the minimum principle ([L]). Thus symplectic 
reduction and the minimum principle are compatible procedures.

\smallskip
For the remainder of this section we assume now that $Z$ is a 
holomorphic $G$-manifold such that there is almost a quotient $Z\qq G$.
More precisely we will assume that $Z\qq G$ is a complex 
space, $\pi:Z\to Z\qq G$ is a surjective $G$-invariant holomorphic
map and there is an analytically  Zariski open $\pi$-saturated
subset  $Z^0$ of $Z$ such that $\pi: Z^0\to Z^0\qq G$ is a geometric 
quotient. Moreover we assume that $\pi: Z^0\to Z^0\qq G$
is a holomorphic fibre bundle. 

\smallskip
Now let $X$ be an invariant domain in $Z$ with $\pi(X)=Z$ and assume that 
$X^0:=X\cap Z^0$ is weakly orbit connected. Let $\phi$ be a $G_\R$-invariant 
plurisubharmonic function such that $\phi^0:=\phi\vert X^0$ is smooth, 
strictly plurisubharmonic on the local $G$-orbits in $X^0$ and an exhaustion mod 
$G_\R$ along $\pi\vert Z^0$. 
Thus the restriction  $\phi^0:=\phi\vert M_\phi^0, \ M_\phi^0:=M_\phi\cap Z^0$ 
induces a plurisubharmonic 
function $\psi^0: Z^0\qq G\to \R$.

\nbbigskip
{\bf Lemma 3.} {\it There is a unique $G$-invariant plurisubharmonic 
function $\Psi:Z\to [-\infty,+\infty)$ which extends 
$\Psi^0:=\psi^0\circ \pi\vert Z^0$.}

\nbigskip
{\it Proof.} The function  $\Psi(z)=\buildunder {g\in \Omega_z}\inf
\phi(g\cdot z)$ is upper semi-continuous on $Z$ where 
$\Omega_z:=\{g\in G;\, g\cdot z\in X\}$. Now 
$\Psi=\Psi^0$ on $Z^0$ (Lemma 2), and $Z\setminus Z^0$ 
is a proper analytic subset of $Z$. Thus $\Psi$ is plurisubharmonic
and by definition $G$-invariant.  
\qed

\nbbigskip
{\it Remark 3.} If $Z\qq G$ is smooth and $\pi$ is an open map, 
then $\psi^0$ extends uniquely to a plurisubharmonic function $\psi$ on 
$Z\qq G$. Of course in this case we have 
$\psi(q)=\buildunder {x \in F_q}\inf
\phi(x)$, where $F_q:=\pi^{-1}(q)\cap X$. If $\phi\vert F_q$ is an exhaustion 
mod $G_\R$, $M_\phi$ intersects every $G_\R$-stable closed analytic subset 
of $F_q$ non trivially. But it might happen that $M_\phi\cap F_q$ is a 
union of several $G_\R$-orbits. On the other hand for $q\in Z^0\qq G$ the 
intersection is exactly one $G_\R$-orbit.

\medskip
Assume now in addition that $Z$ is an open $G$-stable subspace
of a holomorphic Stein $G$-manifold $V$ which is saturated with respect to
$\pi:V\to V\qq G$. We say that $\phi:X\to \R$ is a weak exhaustion of $X$ 
over $V\qq G$ if $\lim \sup \phi(z_n)=+\infty$ for any sequence $(z_n)$ in 
$X$ such that $(\pi(z_n))$ converges to some $q_0$ in the boundary 
$\partial (Z\qq G)$ in $V\qq G$.

\nbbigskip
{\bf Theorem.} {\it Let $Z$ be a $G$-stable $\pi$-saturated 
open subspace of $V$, $X$ an invariant domain in $Z$ with
$G\cdot X=Z$  and $\phi:X\to \R$ a $G_\R$-invariant plurisubharmonic
function. Assume that 
{\parindent=1cm
\smallskip
\litem{{\rm(i)}} $X^0$ is weakly orbit connected, 

\smallskip
\litem{{\rm(ii)}} the restriction of $\phi^0:=\phi\vert X^0$ to the local $G$-orbits 
is strictly plurisubharmonic, 

\smallskip
\litem{{\rm(iii)}} $\phi^0$ is an exhaustion mod $G_\R$ along $\pi\vert Z^0$
and 
\smallskip
\litem{{\rm(iv)}} $\phi$ is a weak exhaustion of $X$ over $V\qq G$,
}
\smallskip\noindent
Then $Z=G\cdot X$ is a Stein manifold.}

\nbigskip
{\it Proof.} Let $z_0\in\partial Z$ and $z_n\in Z$ be such that $z_0=\lim z_n$. 
We have to show $\lim \sup \Psi(z_n)=+\infty$. Thus assume
that $\Psi(z_n)<r$ for all $n$ and some $r\in \R$. There are 
$w_n\in G\cdot M_\phi^0=Z^0$ such that $\Psi(w_n)<r$ and $z_0=\lim w_n$. 
Let $w_n=g_n\cdot x_n$ where $g_n\in G$ and $x_n\in M_\phi^0$. Now 
$\Psi(w_n)=\Psi(x_n)=\phi(x_n)<r$ and, since $Z=G\cdot X$ is saturated, 
$\pi(x_n)=\pi(w_n)\to \pi(z_0)\in \partial (Z\qq G)$. This contradicts 
the assumption that $\phi$ is a weak exhaustion. Thus $Z$ is a domain 
in a Stein manifold with a plurisubharmonic weak exhaustion 
function and therefore Stein.
\qed

\bbigskip
Now let $G$ be complex reductive group and assume that
the semistable quotient $\pi:Z\to Z\qq G$ exists (see [H-M-P]). Thus 
$Z\qq G$ is a complex space whose structure sheaf ${\cal O}_{Z\qq G}(U)
={\cal O}_Z(\pi^{-1}(U)^G$ is the sheaf of invariants and every point in
$Z\qq G$ has an open Stein neighborhood such that the inverse image in $Z$ 
is Stein. 
For example, if $V$ is a holomorphic Stein $G$-manifold, then a 
semistable quotient $V\qq G$ alway exists. Moreover it is shown in [H-M-P]
that $Z$ is a Stein space if and only if $Z\qq G$ is a Stein space.

\smallskip
Assume that $Z$ is connected and that some orbit of maximal 
dimension is closed. Then there exists a proper analytic subset 
$A$ in $Z\qq G$ such that  $Z^o\qq G=Z\qq G\setminus A$ is a 
geometric quotient of $Z^o:=\pi^{-1}(Z\qq G\setminus A)$. In particular,
every fibre of $\pi\vert Z^o$  is $G$-homogeneous or equivalently the 
dimension of the $G$-orbits in $Z^o$ is constant. Every $x\in Z^o$ has a 
$G$-stable neighborhood $U$ which is $G$-equivariantly biholomorphic to 
$G\times_H S$ where $H$ is the isotropy group of $G$ at $x$ and $S$ is a 
Stein space such that the connected component $H^0$ of the identity of 
$H$ acts trivially on $S$. Here $G \times_H S$ denotes the bundle associated
to the $H$-principal bundle $G\to G/H$. Thus locally
$Z^o\qq G$ is given by $S/\Gamma$ where $\Gamma:=H/H_0$ is a 
finite group. Moreover, there is an analytically Zariski open
$G$-stable subset $Z^{oo}$ of $Z$ which is contained in $Z^o$ such 
that the isotropy type is constant. This implies that $Z^{oo}$ is a 
fibre bundle over $Z^{oo}\qq G\subset Z\qq G$.

\nbbbigskip
{\bf 3. Orbit geometry of the future tube.}

\nbbigskip
In the following it will be convenient to introduce a
linear coordinate
change such that $<z,z>=(z_0)^2-(z_1)^2-(z_2)^2-(z_3)^2$
has the form $z_0z_1-z_2z_3$. Thus we set
$$Z:=\biggl({x\atop z}{y\atop w}\biggr)=
\biggl({z_0+z_3\atop z_1+iz_2}{z_1-iz_2\atop z_0-z_3}\biggr)\,$$
and obtain $\det Z=<z,z>$ and $\det \im Z=<\im z,\im z>$ where 
$\im Z:={1\over 2i}(Z - {\bar Z^t})$. 

\smallskip
Let $\H:=\{Z\in V;\ \im Z >0\}$ denote
the generalised upper half plane where $V:=\C^{2\times 2}$.
Note that $\H$ is just the tube over the positive light cone in the new 
coordinates. Moreover $\H$ is stable with respect to the action of 
$G_\R:=\Sl_2(\C)$ which is
given by $G_\R\times \H\to \H,\ (g,Z)\to g*Z:=gZ\bar g^t$. This action is
not effective. The ineffectivity consists of $\Gamma=\{+I, -I\}$ and the 
quotient $\Sl_2(\C)/\Gamma$ is the connected component of the identity of
the homogeneous Lorentz group. 

\smallskip
Let $\H^N:=\H\times\cdots\times \H\subset V\times\cdots\times V=:V^N$ 
denote the $N$-fold product of $\H$ and set 
$G:=(G_\R)^\C=\Sl_2(\C)\times\Sl_2(\C)$
where $G_\R$ is embedded in $G$ via $g\to (g,\bar g)$. The 
diagonal $G_\R$ action on $V^N$ extends to a holomorphic $G$ action
$G\times V^N\to V^N,\ ((g,h), Z^1,\ldots, Z^N)\to (g,h)*(Z^1,\ldots, Z^N):=
(gZ^1h^t,\ldots, gZ^Nh^t)$.

\nbbigskip
{\bf Theorem.} {\it The extended future tube $(\H^N)^\C:=G*\H^N$ is a
domain of holomorphy.}

\nbbigskip
In the proof we will make an axiomatic use of the following statements

\nmedskip
{\bf Fact 1} (see Streater Wightman [S-W], p. 66). The set $\H^N$ is orbit connected 
in $V^N$, i.e., $\{g\in G;\ g*Z\in \H^N\}$ is connected for every $Z\in V^N$.

\nmedskip
{\bf Fact 2.} The extended future tube $G*\H^N$ is saturated with respect 
to $\pi:V^N\to V^N\qq G$.

\nmedskip
Fact 2 implies that the semistable quotient $G*\H^N\qq G$ 
exists and is an open subset of $V^N\qq G$. The quotient map is given
by restricting $\pi:V^N\to V^N\qq G$ to $G*\H^N$. 

\smallskip
There does not seem to be a proof in the literature of Fact 2 but there is
a detailed proof for the whole complex orthogonal group in [H-W]. A slight 
modification of the proof there  can be used for a proof of Fact 2.
In order to be complete let us recall briefly the main steps.
First we note that it is sufficient to show the following (see e.g. [H]).

\nmedskip
{\it Claim.}
If $Z\in \H^N$, then the unique closed orbit $G*W$ 
in the closure of $\overline {G*Z}$ lies in $G*\H^N$.

\nmedskip
This can be seen as follows.
Let $<\ ,\ >$ be the complex Lorenz product, i.e., 
the symmetric bilinear form on $V$ which is
associated to the quadratic form $\det:V\to V$. Thus $V$ is just the standard 
representation of $\tilde G:={\rm O}_4(\C)$. Note that $\tilde G$ has two 
connected components and the connected component of the identity is $G$.
The functions $(Z^1,\ldots, Z^N) \to <Z^i,Z^j>$,
form a set of generators for the algebra of the
$\tilde G$-invariant polynomials on $V^N$. 
Thus the image of $V^N$ in the set of symmetric $N\times N$-matrices of the 
map $\tilde\pi$ which sends $(Z^1,\ldots, Z^N)$ to the matrix  
$(<Z^i,Z^j>)$ is an 
affine variety which is isomorphic to $V^N\qq \tilde G$. 

\smallskip
The matrices 
of rank $3$ or $4$ correspond to fibres of $\tilde\pi$ which are closed 
$\tilde G$-orbits. It follows that the $G$-orbit through every point
$Z\in \H^N$ such that the rank $r$ of $\tilde Z$ is greater or equal to 3
is already closed. Now assume that $r\le 2$. In this case
the following is shown in [H-W]:
There exists an $g\in \tilde G$, $\alpha_j\in \C$ and an $\omega\in V$
with $<\omega, \omega>=0=<\omega, W^j>$ 
such that 
$$Z^j=g*W^j+\alpha_j \omega\ ,\ j=1,\ldots, N\,.$$ 
The proof actually shows that one can choose $g\in G$, i.e., $\det g=1$.
Now an argument of Hall-Wightman ([H-W], p.21) implies that 
$g*W^j\in\H$ for all $j$, i.e., $G*W\subset G*\H^N$.

\nmedskip
{\bf Fact 3.} The function $\phi:\H^N\to \R,\ \phi(Z^1,\ldots, Z^N):=
{1\over \det\, \im Z^1}+\cdots+{1\over \det\, \im Z^N}$ 
is $G_\R$-invariant and strictly plurisubharmonic. Moreover, $\phi$ is a 
weak exhaustion of $\H^N$. 

\nmedskip
The simplest way to see that $\phi$ is strictly plurisubharmonic is
to note that $Z^j\to {1\over \det\im Z^j}$ it is given by the 
Bergmann kernel function on $\H$. 
Since $\det\im Z=0$ for $Z\in\partial\H$,  $\phi$ is a weak 
exhaustion of $\H^N$, i.e.,  
$\phi(Z_k)\to +\infty \ \ 
\hbox{if}\ \ \lim Z_k=Z_0\in \partial (\H^N)\subset V^N$.

\medskip
Let $K_\R:=\{(a,\bar a);\ a\in \Su_2(\C)\}$ be the maximal compact subgroup
of $G_\R$. We set $V^0:=\{Z\in V;\, \det Z\ne 0\}$. Note that $V\qq G\cong \C$
and that after this identification the quotient map is given by
$\det:V\to \C$. In particular, $V^0$ is saturated with respect to 
$V\to V\qq G$. The following Lemma 1 and Lemma 2 can be found in [Z].

\nbbigskip
{\bf Lemma 1.} {\it Let $(W_n)$ be a sequence in $\H$ such 
that $(\pi(W_n))$ converges in $V\qq G$. Then there exist $h_n\in G_\R$
such that a subsequence of $(h_n*W_n)$ converges in $V$.}

\nbigskip
{\it Proof.} There exist $u_n\in K_\R$ such that 
$$X_n:=u_n*W_n=:\bigl({x_n\atop 0}{z_n\atop y_n}\bigr)\,.$$
Since $(\pi(W_n))=(\pi(X_n))$ converges, it follows that 
$\vert \det X_n\vert=\vert x_n y_n\vert\le R$ for some 
$R\ge 0$ and all $n$. Furthermore, $X_n\in \H$ implies that 
${1\over 4}\vert z_n\vert^2 < 
\im x_n \im y_n \le \vert x_n  y_n\vert =\vert \det X_n\vert$.
Therefore $(z_n)$ is bounded. Now $0<\vert x_n y_n\vert\le R$ implies that
$\vert r_n^2 x_n\vert=\vert r_n^{-2} y_n\vert$ for some $r_n>0$. In 
particular the sequence $(r_n^2 x_n, r_n^{-2} y_n)$ is bounded.
Hence $h_n*W_n$ has a convergent subsequence where 
$h_n:=r_n\cdot u_n\in G_\R$ and $r_n$ is identified with 
$\biggl(\bigl({r_n\atop 0}{0\atop {1\over r_n}}\bigr),
\bigl({r_n\atop 0}{0\atop {1\over r_n}}\bigr)\biggr)$.  
\qed

\nbbigskip
{\it Remark.} Geometrically Lemma 1 asserts that $\H$ is relatively compact
over $V\qq G$ mod $G_\R$.

\nbbigskip
{\bf Lemma 2.} {\it Let $(Z_n, W_n)$ be a sequence of points
in $\H\times \H$ and assume that 
{\parindent=1cm
\smallskip
\litem{{\rm (i)}} $\pi(Z_n,W_n)$ converges in $(V\times V)\qq G$ and

\smallskip
\litem{{\rm (ii)}} $W_0=\lim W_n$ exists in $\H$.
}
\smallskip\noindent
Then a subsequence of $(Z_n)$ converges to a $Z_0\in \overline \H$.}

\nbigskip
{\it Proof.} Note that $V\times V^0$ is an open $G$-stable subset 
of $V\times V$ which is saturated with respect to 
$V\times V\to (V\times V)\qq G$ and contains $\H\times \H$. The map 
$V\times V^0\to V,\ (Z,W)\to ZW^{-1}$,
is $G$-equivariant, where $G$ acts on the image $V$ by conjugation 
with the first component, i.e. by $\int(g,h)\cdot X=gXg^{-1}$. 
It is sufficient to show the following 

\nmedskip
{\it Claim.} A subsequence of $(X_n)$ converges.

\nmedskip
Since the image of $X_n:=Z_nW_n^{-1}$ in $V\qq \int\, G$ converges,
the trace and the determinant of $X_n$ and therefore the eigenvalues 
of $X_n$ are bounded. Let 
$u_n=(a_n,\bar a_n)\in K_\R$
be such that $\int\, a_n\cdot X_n=(u_n*Z_n)(u_n*W_n)^{-1}=
\bigl({x_n\atop 0}{z_n\atop y_n}\bigr)$. Since $K_\R$ is compact,
we may assume that $X_n=\bigl({x_n\atop 0}{z_n\atop y_n}\bigr)$.

\smallskip
Let $W_n=:\bigl({a_n\atop c_n}{b_n\atop d_n}\bigr)$ and 
$W_0=:\bigl({a_0\atop c_0}{b_0\atop d_0}\bigr)$. By assumption we have
$W_0 \in \H$. Therefore $\im d_0\ne 0$. From 
$$Z_n=X_nW_n=
\left(\matrix{
x_n a_n+z_n c_n & x_n b_n+z_n d_n \cr
y_n c_n & y_n d_n \cr }\right)\in \H$$ 
it follows that 
$${1\over \vert z_n\vert^2}\bigl(
\im(x_n a_n+ z_n c_n)\im (y_n d_n)-{1\over 4}\vert x_n b_n + z_n d_n
-\bar y_n \bar c_n\vert^2\bigr)>0$$
for $z_n\ne 0$.
Since the eigenvalues $x_n,y_n$ and $a_n, b_n, c_n, d_n$ are bounded, 
$d_0\ne0$ implies that $\vert z_n\vert$ is bounded. Thus 
$(X_n)$ has a convergent subsequence.
\qed

\nbbigskip
In the above proof we used that $\H\subset V^0$
which is implied by $\det\im Z\le \vert \det Z\vert$.

\nbbigskip
{\bf Corollary 1.} {\it If $Z_n=(Z_n^1,\dots,Z_n^N)\in \H^N$ are such that
$(\pi(Z_n))$ converges in $V^N\qq G$ and  
$(Z_n^N)$ converges in $\H$, then  $(Z_n)$ has a convergent subsequence in 
$\overline{\H}^N $.}  
\qed

\nbbigskip
{\bf Lemma 3.} {\it $\phi$ is a weak exhaustion of $X$ over $V\qq G$.}

\nbigskip
{\it Proof.} Let $(Z_n)=((Z_n^1,\ldots ,Z_n^N))$ be a sequence in $\H^N$ 
such that 
$q:=\lim \pi(Z_n)\in \partial (G*\H^N\qq G)\subset V^N\qq G$ exists. 
There are $h_n\in G_\R$ such that a subsequence of $(h*Z^N_n)$ 
converges to $W^N \in \overline \H$ (Lemma 1). Now, if $W^N\in \partial \H$,
then $\lim \sup \phi(Z_n)=+\infty$. Thus assume that $W^N\in \H$. It follows
that $(h_n*Z_n)$ has a subsequence which converges to $W\in \overline \H^N$
(Corollary 1). But $W$ is not in $\H^N$, since 
$q=\pi(W) \in \partial (G*\H^N\qq G)$. Thus $W\in \partial \H^N$ and therefore
again $\lim \sup \phi(Z_n)=+\infty$ follows.
\qed

\nbbigskip
{\bf Lemma 4.} {\it The function $\phi$ is an exhaustion mod $G_\R$ along 
$\pi$.}

\nbigskip
{\it Proof.} For $r>0$ let $Z_n\in \H^N$, $Z_n=:(Z_n^1,\ldots,Z_n^N)$,
be such that $\phi(Z_n)\le r$ and assume that $\lim \pi(Z_n)$ exists in 
$G*\H^N\qq G$. Thus there are $h_n\in G_\R$ such that $(h_n*Z_n^N)$ has 
a subsequence which converges to some $W^N\in \overline \H$. 
If $W^N\in\partial \H$, then $\phi(Z_n)$ goes to infinity. This
contradicts $\phi(h_n*Z_n)\le r$. Thus $W^N\in \H$ and therefore 
$(h_n*Z_n)$ has a subsequence with limit 
$W=(W^1,\ldots, W^N)\in \overline \H^N$. The same argument as above implies 
that $W^j\in \H$ for $j=1,\ldots, N$.
\qed

\nbbigskip
{\it Proof of the Theorem.} 
From the invariant theoretical point of view the 
$G=\Sl_2(\C)\times \Sl_2(\C)$ action
on $V^N$ is the $N$-fold product of the standard representation 
of $\hbox{SO}_4(\C)$ on $\C^4$. It is well known that for any $N=1,2,\ldots$ 
the generic $G$-orbit in $V^N$ is closed. Let $(V^N)^0$ denote the set of
points in $V^N$ which lie in a generic closed orbit, i.e., 
$(V^N)^0$ is a  union of the fibres of the
quotient $V\to V\qq G$ which consist exactly of one $G$-orbit. 
Since the $G_\R$-action on $\H$ is proper, $G_\R$ acts properly 
on $\H^N$. It follows from the results in \S2 that there is a $G$-invariant  
plurisubharmonic function $\Psi$ on $G*\H^N$ which is a weak exhaustion.
Thus $G*\H^N$ is a domain of holomorphy.
\qed

\nbbigskip
{\bf Corollary 2.} {\it The image $G*\H^N\qq G$ of $\H^N$ in 
$V^N\qq G$ is an open Stein subspace.}
\qed

\nbbbigskip
{\it References}

{\parindent=1,5cm

\bigskip
\litem{[A]} Ammon, M.: 
{\it Komplexe Strukturen auf Quotienten von Kempf-Ness Mengen.}  
Dissertation Bochum Januar 1997

\smallskip
\litem{[A-H-H]} Ammon, M., Heinzner, P.; Huckleberry, A. T.: 
{\it Complex structures on symplectic reduction spaces.}  
In preparation

\smallskip
\litem{[B-L-T]} Bogolyubov, N. N., Logunov, A. A.; Todorov, I. T.: 
Introduction to axiomatic quantum field theory. Benjamin, Reading, Mass.,
1974

\smallskip
\litem{[H-W]} Hall, D.; Wightman, A. S.: {\it A theorem on invariant analytic 
functions with applications to relativistic quantum field theory.}
Kgl. Danske Videnskab. Selskab, Mat. Fys. Medd. 31, Nr. 5 (1957), 1--41

\smallskip
\litem{[H]} Heinzner, P.: {\it Geometric invariant theory on Stein spaces.}
Math. Ann. 289 (1991), 631-662

\smallskip
\litem{[H-H-K]} Heinzner, P.; Huckleberry, A. T.; Kutzschebauch, F.: {\it 
Abels' Theorem in the real analytic case and applications to 
complexifications.} In: Complex Analysis and Geometry, 
Lecture Notes in Pure and Applied Mathematics, Marcel Decker 1995, 229--273

\smallskip
\litem{[H-H-L]} Heinzner, P.; Huckleberry, A. T.; Loose, F.: {\it K\"ahlerian
extensions of the symplectic reduction.}  
J. reine und angew. Math. 455,
123--140 (1994)

\smallskip
\litem{[H-L]} Heinzner, P.; Loose, F.: {\it A global slice for 
proper Hamiltonian actions.} 
Preprint Bochum 1997, 8p.

\smallskip
\litem{[H-M-P]} Heinzner, P.; Migliorini, L.; Polito, M.: {\it Semistable quotients.}  
Annali della Scuola Normale Superiore di Pisa (1997) (to appear) 16p.

\smallskip
\litem{[H-S]} Heinzner, P.; Sergeev, A.: {\it The extended matrix disc is
a domain of holomorphy.} 
Math. USSR Izvestiya 38, 637--645 (1992)

\smallskip
\litem{[J]} Jost, R.: The general theory of quantized fields. Lecture Notes in 
Appl. Math. Vol. 4 Amer. Math. Soc., Providence, R.I.,  1965

\smallskip
\litem{[L]} Loeb, J.-J.: {\it Action d'une forme r\'eelle d'un groupe de
Lie complexe sur les fonctions plurisousharmoniques.} 
Ann. Inst. Fourier, Grenoble 35, 4, 59--97 (1985)

\smallskip
\litem{[S-V]} Sergeev, A. G.; Vladimirov, V. S.: {\it  Complex analysis in 
the future tube.} In: Encyclopedia of Math. Sci. vol. 8 (Several Complex
Variables II), Springer Verlag

\smallskip
\litem{[S-L]} Sjamaar, R.; Lerman, E.: {\it Stratified symplectic
	spaces and reduction.} Annals  of  Math. 134, 375--422 (1991)

\smallskip
\litem{[S-W]} Streater, R.F.; Wightman, A.S.: PCT spin statistic, and all that.
W. A. Benjamin, INC., New York Amsterdam 1964

\smallskip
\litem{[Z]} Zhou, Xiang-Yu: {\it On the extended future tube conjecture.}
Preprint and supplementary pages 1997

\par}

\bigskip
Peter Heinzner

Brandeis University

Department of Mathematics

Waltham, MA 02254-9110

USA

\end